\documentclass[10pt,reqno]{amsart}
\usepackage{amsmath}
\usepackage{amsthm}
\usepackage{amssymb}
\usepackage{graphicx}
\usepackage{amsfonts}
\usepackage{latexsym}

\newcommand{\dist}{\operatorname{dist}}

\newcommand{\B}{ \mathcal{B}}
\newcommand{\C}{ \mathbb{C}}
\newcommand{\D}{ \mathbb{D}}
\newcommand{\dD}{ \partial\mathbb{D}}

\newcommand{\T}{ \mathcal{T}}

\newcommand{\ran}{\operatorname{ran}}

\newcommand{\norm}[1]{\left\| #1 \right\|}
\newcommand{\inner}[1]{\left< #1 \right>}

\newcommand{\N}{\mathbb{N}}

\newcommand{\h}{\mathcal{H}}
\newcommand{\K}{\mathcal{K}}

\renewcommand{\phi}{\varphi}
\newcommand{\Aut}{\operatorname{Aut}}

\numberwithin{equation}{section}
\theoremstyle{plain}
\newtheorem{Proposition}[equation]{Proposition}
\newtheorem{Corollary}[equation]{Corollary}
\newtheorem{Theorem}[equation]{Theorem}
\newtheorem{Lemma}[equation]{Lemma}
\theoremstyle{definition}
\newtheorem{Definition}[equation]{Definition}

\newtheorem{Remark}[equation]{Remark}
\newtheorem{Question}[equation]{Question}

\allowdisplaybreaks

\begin{document}

\bibliographystyle{amsplain}

\title[Truncated Toeplitz operators]
{Truncated Toeplitz Operators:\\Spatial Isomorphism, Unitary Equivalence,\\and Similarity}

\author[J.A.~Cima]{Joseph A. Cima}
	\address{Department of Mathematics, University of North Carolina, Chapel Hill, North Carolina 27599}
	\email{cima@email.unc.edu}

\author[S.R.~Garcia]{Stephan Ramon Garcia}
	\address{Department of Mathematics, Pomona College, Claremont, California 91711}
	\email{Stephan.Garcia@pomona.edu}

\author[W.T.~Ross]{William T. Ross}
	\address{Department of Mathematics and Computer Science, University of Richmond, Richmond, Virginia 23173}
	\email{wross@richmond.edu}

\author[W.R.~Wogen]{Warren R. Wogen}
	\address{Department of Mathematics, University of North Carolina, Chapel Hill, North Carolina 27599}
	\email{wrw@email.unc.edu}
\maketitle

\begin{abstract}
	A \emph{truncated Toeplitz operator}
	 $A_{\phi}:\K_{\Theta} \to \K_{\Theta}$ is the compression of a
	Toeplitz operator $T_{\phi}:H^2\to H^2$ to a model space
	$\K_{\Theta} := H^2 \ominus \Theta H^2$.
	For $\Theta$ inner, let $\T_{\Theta}$ denote the set of all
	bounded truncated Toeplitz operators on $\K_{\Theta}$.
	Our main result is a necessary and sufficient condition
	on inner functions $\Theta_1$ and $\Theta_2$ which guarantees that
	$\mathcal{T}_{\Theta_1}$ and $\mathcal{T}_{\Theta_2}$ are spatially isomorphic.
	(i.e., 	$U\T_{\Theta_1} = \T_{\Theta_2}U$ for some unitary $U:\K_{\Theta_1} \to \K_{\Theta_2}$).
	We also study operators which are unitarily equivalent to truncated Toeplitz operators
	and we prove that every operator on a finite dimensional Hilbert
	space is similar to a truncated Toeplitz operator.
\end{abstract}

\section{Introduction}
In this paper we consider several questions concerning spatial isomorphism, unitary equivalence,
	and similarity in the setting of truncated Toeplitz operators.  Loosely put, a \emph{truncated Toeplitz operator}
	is the compression $A_{\phi}:\K_{\Theta} \to \K_{\Theta}$ of a standard
	Toeplitz operator $T_{\phi}:H^2\to H^2$ to a Jordan model space
	$\K_{\Theta} := H^2 \ominus \Theta H^2$ (here $\Theta$ denotes an inner function).
	We discuss these definitions and the related preliminaries in Section \ref{SectionPrelim}.
	The reader is directed to the recent survey of Sarason \cite{Sarason}
	for a more thorough account.
	
	For a given inner function $\Theta$, we let $\T_{\Theta}$ denote the set of all
	bounded truncated Toeplitz operators on $\K_{\Theta}$.
	The main result of the paper (Theorem \ref{Main-SI}) is a simple necessary and sufficient condition
	on inner functions $\Theta_1$ and $\Theta_2$ which guarantees that the corresponding spaces
	$\mathcal{T}_{\Theta_1}$ and $\mathcal{T}_{\Theta_2}$ are spatially isomorphic (i.e.,
	$U\T_{\Theta_1} = \T_{\Theta_2}U$ for some unitary $U:\K_{\Theta_1} \to \K_{\Theta_2}$).
	This result and its ramifications are discussed in Section \ref{SectionMain} while the
	proof is presented in Section \ref{SectionProof}.

	In Section \ref{SectionUE}, we study the operators which
	are unitarily equivalent to truncated Toeplitz operators (UETTO).
	The class of such operators is surprisingly large and includes, for instance, the
	Volterra integration operator \cite{MR0192355}.
	We add to this class by showing that several familiar classes of
	operators (e.g., normal operators) are UETTO.

	We conclude this paper in Section \ref{SectionSimilar}
	by showing that every operator on a finite dimensional Hilbert
	space is similar to a truncated Toeplitz operator (Theorem \ref{Jordan-TTO}).
	In other words, we prove that the inverse Jordan structure problem is always solvable in
	the class of truncated Toeplitz operators.  This stands in  contrast to the situation for
	Toeplitz \emph{matrices} \cite{MR1839449}.

\section{Preliminaries}\label{SectionPrelim}

In the following, $H^2$ denotes the classical
Hardy space on the open unit disk $\D$ \cite{Duren, MR2261424}.
The unit circle $|z|=1$ is denoted by $\dD$ and we let
$L^2 := L^2(\dD)$ and $L^{\infty} := L^{\infty}(\dD)$ denote the
usual Lebesgue spaces on $\dD$.

\subsection*{Model spaces}
	To each non-constant inner function $\Theta$ there corresponds
	a \emph{model space} $\mathcal{K}_{\Theta}$ defined by
	\begin{equation}\label{eq-Alternate}
		\mathcal{K}_{\Theta} := H^2 \ominus \Theta H^2.
	\end{equation}
	This terminology
	stems from the important role that $\mathcal{K}_{\Theta}$ plays in the
	\emph{model theory} for Hilbert space contractions -- see \cite[Part C]{MR1892647}.

	The kernel functions
	\begin{equation}\label{eq-ReproducingKernel}
		K_{\lambda}(z)
		= \frac{1 - \overline{\Theta(\lambda)} \Theta(z)}{1 - \overline{\lambda} z}, \quad z,\lambda \in \D,
	\end{equation}
	belong to $\mathcal{K}_{\Theta}$ and
	enjoy the \emph{reproducing property}
	\begin{equation}\label{eq-ReproducingProperty}
		\inner{f, K_\lambda} = f(\lambda), \quad \lambda \in \D, \,f\in \mathcal{K}_{\Theta}.
	\end{equation}
	If $\Theta$ has an angular derivative in the sense of Carath\'eodory
	(ADC) at $\lambda \in \dD$ \cite[Sect.~2.2]{Sarason}
	then $K_{\lambda}$ belongs to $\mathcal{K}_{\Theta}$ and the formulae \eqref{eq-ReproducingKernel}
	and \eqref{eq-ReproducingProperty} still hold.
	Letting $P_{\Theta}$ denote the orthogonal projection of $L^2$ onto $\mathcal{K}_{\Theta}$,
	we observe that
	\begin{equation} \label{eq-ProjectionTheta}
		[P_{\Theta} f](\lambda) = \inner{f, K_{\lambda}}, \quad f \in L^2, \lambda \in \D.
	\end{equation}
	 The preceding formula
	remains valid for $\lambda \in \dD$ so long as $\Theta$ has an ADC there.

	We let
	\begin{equation} \label{nrk}
		k_{\lambda} := \frac{K_{\lambda}}{\norm{K_{\lambda}}}
	\end{equation}
	denote the \emph{normalized reproducing kernel} at $\lambda$ and, when we wish to be specific
	about the underlying inner function $\Theta$ involved, we write
	$K_{\lambda}^{\Theta}$ and $k_{\lambda}^{\Theta}$ in place of
	$K_{\lambda}$ and $k_{\lambda}$, respectively.	
	
	There is a natural \emph{conjugation} (a conjugate-linear isometric involution) on $\K_{\Theta}$
	defined in terms of boundary functions by
	\begin{equation}\label{eq-ModelConjugation}
		C f := \overline{ f z} \Theta.
	\end{equation}
	Although at first glance the expression $\overline{fz}\Theta$ in \eqref{eq-ModelConjugation}
	does not appear to correspond to the boundary values of an $H^2$ function, let alone
	one in $\mathcal{K}_{\Theta}$, a short computation using \eqref{eq-Alternate}
	reveals that if $f \in \mathcal{K}_{\Theta}$ and $h \in H^2$, then
	$\inner{Cf, \Theta h} = 0 = \inner{C f, \overline{z h}}$
	whence $Cf$ indeed belongs to $\mathcal{K}_{\Theta}$.
		
	A short calculation reveals that
	\begin{equation*}
		[C K_{\lambda}](z) = \frac{\Theta(z) - \Theta(\lambda)}{z - \lambda}.
	\end{equation*}
	Moreover, the preceding also holds for $\lambda \in \dD$ so long as $\Theta$ has an ADC there.

\subsection*{Truncated Toeplitz operators}

	Since $\K_{\Theta}$ is the closed linear span of the backward shifts
	$S^{*} \Theta, S^{* 2} \Theta, \ldots$ of $\Theta$ \cite[p.~83]{CR}, where $S^{*} f = (f - f(0))/z$, it follows that the subspace
	\begin{equation*}
		\K_{\Theta}^{\infty} := \mathcal{K}_{\Theta} \cap H^{\infty}
	\end{equation*}
	of all bounded functions in $\K_{\Theta}$ is dense in $\mathcal{K}_{\Theta}$.


	Keeping these results in mind, for a fixed inner function $\Theta$ and
	any $\phi \in L^2$, the corresponding \emph{truncated Toeplitz operator}
	$A_{\phi}:\K_{\Theta} \to \K_{\Theta}$ is the densely defined operator
	\begin{equation} \label{eq-TTODefinition}
		A_{\varphi} f = P_{\Theta}(\varphi f).
	\end{equation}
	When we wish to be specific about the underlying inner function $\Theta$, we use
	the notation $A^{\Theta}_{\phi}$ to denote the truncated
	Toeplitz operator with symbol $\phi$ acting on the model space $\mathcal{K}_{\Theta}$.
	In most cases, however, $\Theta$ is clear from context and we simply write $A_{\phi}$.
	
	Although one can pursue the subject of unbounded truncated Toeplitz operators
	much further \cite{MR2468883, MR2418122}, we are concerned here
	with those which have a bounded extension to $\K_{\Theta}$.
	
	\begin{Definition}
		Let $\mathcal{T}_{\Theta}$ denote the set of all
		truncated Toeplitz operators which extend boundedly to all of $\mathcal{K}_{\Theta}$.
	\end{Definition}
	
Certainly $A_{\phi} \in \mathcal{T}_{\Theta}$ when $\phi \in L^{\infty}$. However \cite[Thm.~3.1]{Sarason}, there are an abundance of unbounded $\phi \in L^{2}$ for which $A_{\phi} \in \mathcal{T}_{\Theta}$.
It is important to note that $\T_{\Theta}$ is not an algebra since the product of truncated
	Toeplitz operators need not be a truncated Toeplitz operator (a simple counterexample
	can easily be deduced from \cite[Thm.~5.1]{Sarason}).  On the other hand, it turns out
	that $\mathcal{T}_{\Theta}$ is a weakly closed linear subspace of the bounded operators on $\mathcal{K}_{\Theta}$ \cite[Thm.~4.2]{Sarason}.
	Moreover, if $\Theta$ is a finite Blaschke product of order $n$, then one can show that
	$\dim \mathcal{T}_{\Theta} = 2 n - 1$ (see Lemma \ref{Sarason-tensors-kck} below).

\subsection*{Complex symmetric operators}
	Of particular importance to the study of truncated Toeplitz operators
	is the notion of a complex symmetric operator \cite{G-P, G-P-II}.  Let us briefly
	discuss the necessary preliminaries.
In the following, we let $\h$ denote
	a separable complex Hilbert space and $\mathcal{B}(\mathcal{H})$ denote the bounded linear operators on $\mathcal{H}$.

	\begin{Definition}
		A \emph{conjugation} on $\h$ is a conjugate-linear operator $C:\h \rightarrow \h$,
		which is both \textit{involutive} (i.e., $C^2 = I$)
		and \textit{isometric} (i.e., $\inner{Cx,Cy} = \inner{y,x}$ for all $x,y \in \h).$
	\end{Definition}

	The standard example of a conjugation is entry-by-entry complex conjugation on an $l^2$-space.
	In fact, each conjugation is unitarily equivalent to the canonical conjugation on a $l^2$-space
	of the appropriate dimension \cite[Lem.~1]{G-P}.
Having discussed conjugations, we next consider certain operators which are compatible
	with them.

	\begin{Definition} \label{def-CSO}
		We say that $T \in \B(\h)$ is \emph{$C$-symmetric}
		if $T^* = CTC$ for some conjugation $C$ on $\h$.
		We say that $T$ is \emph{complex symmetric} if there exists
		a conjugation $C$ with respect to which $T$ is $C$-symmetric.
	\end{Definition}

Recall the conjugation $C$ defined on $\mathcal{K}_{\Theta}$ from \eqref{eq-ModelConjugation}. The following result is from \cite{G-P}.

\begin{Proposition} \label{TTOareCSO}
Every $A_{\phi} \in \mathcal{T}_{\Theta}$ is $C$-symmetric.
\end{Proposition}


\subsection*{Clark operators}
	Let us now review a few necessary facts about the theory of Clark
	unitary operators \cite{MR0301534}.   For a more complete account of this theory we refer the reader to
	\cite{MR2215991, MR2198367, MR1289670}. To avoid needless technicalities,
	we assume that the underlying inner function $\Theta$ satisfies $\Theta(0) = 0$.
	For $\alpha \in \dD$, the operator $U_{\alpha}: \mathcal{K}_{\Theta} \to \mathcal{K}_{\Theta}$
	defined by the formula
	\begin{equation} \label{Clark-U-defn}
		U_{\alpha}f = A_{z}f+ \alpha \inner{f, \overline{z}\Theta} 1, \quad f \in \K_{\Theta},
	\end{equation}
	is called a \emph{Clark operator}.
	One can show that each Clark operator $U_{\alpha}$ is unitary and that
	every unitary rank-one perturbation of the truncated shift operator $A_z$ takes the form
	$U_{\alpha}$ for some $\alpha \in \dD$.  Less well-known is the fact that
	each Clark operator $U_{\alpha}$ on $\K_{\Theta}$
	belongs to $\mathcal{T}_{\Theta}$ \cite[p.~524]{Sarason}.
	There is also the following theorem \cite[p.~515]{Sarason}.

	\begin{Theorem}[Sarason]\label{T-Sarason-commutes}
		If $A$ is a bounded operator
		on $\mathcal{K}_{\Theta}$ which commutes with $U_{\alpha}$
		for some $\alpha \in \dD$, then $A \in \mathcal{T}_{\Theta}$.
	\end{Theorem}

	Since $U_{\alpha}$ is a cyclic unitary operator \cite[Thm.~8.9.10]{MR2215991},
	the Spectral Theorem asserts that there is a measure $\mu_{\alpha}$ on $\dD$
	such that $U_{\alpha}$ is unitarily equivalent to the operator
	$[M_{\zeta} f](\zeta) = \zeta f(\zeta)$ of multiplication by the independent variable $\zeta$
	on $L^{2}(\mu_{\alpha})$.  Moreover,
	the measure $\mu_{\alpha}$ is carried by the set
	\begin{equation*}
		E_{\alpha} := \left\{ \zeta \in \dD: \lim_{r \to 1^{-}} \Theta(r \zeta) = \alpha \right \}
	\end{equation*}
	and is therefore singular with respect to Lebesgue measure on $\dD$.
	The \emph{Clark measure} $\mu_{\alpha}$ constructed above can also easily be
	obtained using the Herglotz Representation Theorem
	for harmonic functions with positive real part \cite[Ch.~9]{MR2215991}. As a consequence of this, one can use the fact that $\Theta(0) = 0$ to see that $\mu_{\alpha}$ is a probability measure.

	It is important to note that the preceding recipe can essentially be reversed.
	We record this observation here for future reference (see \cite[p.~202]{MR2215991} for details).

	\begin{Proposition} \label{converse-Clark}
		If $\mu$ is a singular probability measure on $\dD$,
		then there is an inner function $\Theta$ with $\Theta(0) = 0$
		such that the Clark measure for $\Theta$ at $\alpha  = 1$ is $\mu$.
		In particular, $\mu$ is the spectral measure for the Clark unitary
		operator $U_1$ on $\mathcal{K}_{\Theta}$.
	\end{Proposition}

	In the finite-dimensional case, the Clark measures $\mu_{\alpha}$ can be computed explicitly.
	If $\Theta$ is a finite Blaschke product of order $n$, then $\dim \K_{\Theta} = n$
	and the set $E_{\alpha}$ consists of the $n$ distinct points $\zeta_1, \zeta_2,\ldots,\zeta_n$
	on $\dD$ for which $\Theta(\zeta_j) = \alpha$.  The corresponding normalized reproducing kernels $k_{z_j}$ satisfy
	$U_{\alpha} k_{\zeta_j} = \zeta_j k_{\zeta_j}$ for $j=1,2,\ldots,n$ and form an
	orthonormal basis for the model space $\K_{\Theta}$.

\subsection*{Rank one operators in $\T_{\Theta}$}
	Let us conclude these preliminaries with a few words concerning
	truncated Toeplitz operators of rank one.  First recall that for each pair $f, g$ of vectors in
	a Hilbert space $\h$ the operator $f \otimes g: \mathcal{H} \to \mathcal{H}$ is defined by setting
	\begin{equation} \label{rank-1-tensor-def}
		(f \otimes g)(h) := \inner{h, g}\! f.
	\end{equation}
	Observe that $f \otimes g$ has a rank one and range $\C f$.
	Moreover, we also have
	$\|f \otimes g\| = \norm{f} \norm{g}$.
	The proof of the next lemma is elementary and is left to the reader.

	\begin{Lemma} \label{L-tensors}
		Let $\h_1,\h_2$ be Hilbert spaces and let $f_1,g_1 \in \h_1$
		and $f_2,g_2 \in \h_2$ be unit vectors.
		 \begin{enumerate}\addtolength{\itemsep}{0.5\baselineskip}
			\item If $U: \mathcal{H}_1 \to \mathcal{H}_2$ is a unitary operator such that
				\begin{equation*}
					U (f_1 \otimes g_1) U^{*} = f_2 \otimes g_2,
				\end{equation*}
				then there exists a $\zeta \in \dD$ such that
				$U f_1 = \zeta f_2$ and $U g_1 = \zeta g_2$. In particular, we have
				$\inner{f_1, g_1}_{\mathcal{H}_1} = \inner{f_2, g_2}_{\mathcal{H}_2}$.

			\item Conversely, if $\inner{f_1, g_1} = \inner{f_2, g_2}$,
				then the operators $f_1 \otimes g_1$ and $f_2 \otimes g_2$ are unitarily equivalent.
		\end{enumerate}
	\end{Lemma}

	The following useful lemma completely characterizes the truncated Toeplitz
	operators of rank one \cite[Thm.~5.1, 7.1]{Sarason}.  We remind the reader that
	$K_{\lambda}$ denotes the reproducing kernel \eqref{eq-ReproducingKernel} for $\K_{\Theta}$ and $C$ denotes the conjugation on $\mathcal{K}_{\Theta}$ from \eqref{eq-ModelConjugation}.

	\begin{Lemma}[Sarason] \label{Sarason-tensors-kck}
		Let $\Theta$ be an inner function.
		 \begin{enumerate}\addtolength{\itemsep}{0.5\baselineskip}
			\item For each $\lambda \in \D$, the operators
				$K_{\lambda} \otimes C K_{\lambda}$ and $C K_{\lambda} \otimes K_{\lambda}$
				belong to $\mathcal{T}_{\Theta}$.
				
			\item If $\eta \in \dD$ and $\Theta$ has a ADC at $\eta$, then
				$K_{\eta} \otimes K_{\eta} \in \mathcal{T}_{\Theta}$.

			\item The only rank-one operators in $\T_{\Theta}$ are
				the nonzero scalar multiples of the operators from (i) and (ii). 	
				
			\item If $\Theta$ is a Blaschke product of order $n$, then \smallskip
				 \begin{enumerate}\addtolength{\itemsep}{0.5\baselineskip}
					\item $\dim \T_{\Theta} = 2n - 1$.
					
					\item If $\lambda_1, \lambda_2, \ldots, \lambda_{2n-1}$ are distinct points of $\D$, then the operators
						\begin{equation*}
							K_{\lambda_i} \otimes C K_{\lambda_i}, \quad 1 \leq i \leq 2 n - 1,
						\end{equation*}
						form a basis for $\T_{\Theta}$.
				\end{enumerate}
		\end{enumerate}
	\end{Lemma}
	
	Elementary complex analysis tells us that the automorphism group
	$\Aut(\D)$ of $\D$ can be explicitly presented as
	\begin{equation*}
		\Aut(\D) = \{\zeta \phi_a: \zeta \in \dD, a \in \D\}
	\end{equation*}
	where $\phi_a$ denotes the M\"obius transformation
	\begin{equation}\label{eq-Mobius}
		\phi_{a}(z) := \frac{z - a}{1 - \overline{a} z}.
	\end{equation}
	For an inner function $\Theta$ with $\Theta \notin \Aut(\D)$ we have the following lemma:

	\begin{Lemma} \label{L-LI}
		Suppose $\Theta$ is inner with $\Theta \not \in \operatorname{Aut}(\D)$.\smallskip
		 \begin{enumerate}\addtolength{\itemsep}{0.5\baselineskip}
			\item If $\lambda_1, \lambda_2 \in \D$, then $K_{\lambda_1}$ is not a scalar multiple of $C K_{\lambda_2}$.
			\item If also $\lambda_1 \not = \lambda_2$, then $K_{\lambda_1}$ is not a scalar multiple of $K_{\lambda_2}$.
			\item If $\lambda \in \D$, then $K_{\lambda} \otimes C K_{\lambda}$ is not self-adjoint.
		\end{enumerate}
	\end{Lemma}

	\begin{proof}
		Statements (i) and (ii)
		are easy computations. For (iii), note that (i) shows that
		$K_{\lambda} \otimes C K_{\lambda}$ and $(K_{\lambda} \otimes
		C K_{\lambda})^{*}  = C K_{\lambda} \otimes K_{\lambda}$ have
		different ranges.
	\end{proof}

\section{When are $\T_{\Theta_1}$ and $\T_{\Theta_2}$ spatially isomorphic?}\label{SectionMain}

	In this section we consider the problem of determining when two spaces
	$\T_{\Theta_1}, \T_{\Theta_2}$ of truncated Toeplitz operators are spatially
	isomorphic.  Let us recall the following definition.

	\begin{Definition}
		For $j = 1, 2$, let $\mathcal{H}_j$ be a Hilbert space and $\mathcal{S}_j$
		be a subspace of $\mathcal{B}(\mathcal{H}_j)$. We say that
		$\mathcal{S}_1$ is \textit{spatially isomorphic} to $\mathcal{S}_2$, written
		$\mathcal{S}_1 \cong \mathcal{S}_2$, if there is a unitary operator
		$U: \mathcal{H}_1 \to \mathcal{H}_2$ so that the map
		\begin{equation*}
			S \mapsto U S U^{*}, \quad S \in \mathcal{S}_1,
		\end{equation*}
		carries $\mathcal{S}_1$ onto $\mathcal{S}_2$.
		In this case we often write $U \mathcal{S}_1 U^{*} = \mathcal{S}_2$.
	\end{Definition}
	
	Let us be more precise about our main problem.  Spatial isomorphisms
	of the spaces $\T_{\Theta}$ give rise to an equivalence relation
	on the collection of all inner functions and we wish to determine the
	structure of the corresponding equivalence classes.

	If $\Theta$ is an inner function and $\psi$ belongs to $\Aut(\D)$, then the
	functions $\psi \circ \Theta$ and $\Theta \circ \psi$ are also inner and hence
	\begin{equation*}
		\mathcal{O}(\Theta) := \{\psi_1 \circ \Theta \circ \psi_2\,:\, \psi_1, \psi_2 \in \Aut(\D)\}
	\end{equation*}
	consists precisely of those inner functions that can be obtained from $\Theta$
	by pre- and post-composition with disk automorphisms.  It turns out that while
	$\Theta_1 \in \mathcal{O}(\Theta_2)$ is a sufficient condition for ensuring that
	$\T_{\Theta_1} \cong \T_{\Theta_2}$, it is not necessary.
	To formulate the correct theorem, we introduce the
	conjugation $f \mapsto f^{\#}$ on $H^2$ by setting
	\begin{equation}\label{eq-Sharp}
		f^{\#}(z) := \overline{f(\overline{z})}
	\end{equation}
	and we note that $\Theta^{\#}$ is inner if and only if $\Theta$ is.  Moreover, note that
	the $\#$ operation naturally extends to a conjugation on all of $L^2$. The main theorem of this section is the following:

	\begin{Theorem} \label{Main-SI}
		For inner functions $\Theta_1$ and $\Theta_2$,
		\begin{equation}\label{eq-Implication}
			\mathcal{T}_{\Theta_1} \cong \mathcal{T}_{\Theta_2}
			\quad\Leftrightarrow\quad
			\Theta_1 \in \mathcal{O}(\Theta_2) \cup \mathcal{O}(\Theta_2^{\#}).
		\end{equation}
	\end{Theorem}

	The proof of the preceding theorem is somewhat long and it requires
	a number of technical lemmas.  We therefore defer the proof until Section \ref{SectionProof}.
	
	\medskip

	It is natural to ask if there are simple geometric conditions
	on the zeros of Blashke products $B_1$ and $B_2$ that will ensure that
	$\mathcal{T}_{B_{1}} \cong \mathcal{T}_{B_{2}}$.
	While the general question appears difficult, several partial results
	are available.  For instance, if $B_1$ and $B_2$ are Blaschke products of order $2$, then
	$\mathcal{T}_{B_1} \cong \mathcal{T}_{B_2}$ (see Theorem \ref{T-2x2} below).
	Another special case is handled by the following corollary.  Before
	presenting it, we require a few words concerning the hyperbolic metric on $\D$.
	
	The hyperbolic (or Poincar\'{e}) metric on $\D$ is defined for $z_1,z_2 \in \D$ by
	\begin{equation} \label{hyperbolic}
		\rho(z_1, z_2) = \inf_{\gamma} \int_{\gamma} \frac{2 |d z|}{1 - |z|^2},
	\end{equation}
	where the infimum is taken over all arcs $\gamma$ in $\D$ connecting $z_1$ and $z_2$.
	It is well-known that the hyperbolic metric $\rho$ is conformally invariant in the sense that
	\begin{equation*}
		\rho(z_1, z_2) = \rho(\psi(z_1), \psi(z_2)), \quad \forall \psi \in \Aut(\D).
	\end{equation*}
	Moreover,
	\begin{equation} \label{log-identity}
		\rho(0, z) = \log \frac{1 + |z|}{1 - |z|}
	\end{equation} and the geodesic through $0, z$ turns out to be $[0, z]$,
	the line segment from $0$ to $z$. The reader can
	consult \cite[p.~4]{MR2261424} for further details.

	\begin{Corollary}
		For a finite Blaschke product $B$ of order $n$, we have
		$\mathcal{T}_{z^n} \cong \mathcal{T}_B$ if and only if either
		$B$ has one zero of order $n$ or $B$ has $n$ distinct zeros all lying
		on a circle $\Gamma$ in $\mathbb{D}$ with the property that if these zeros
		are ordered according to increasing argument on $\Gamma$, then
		adjacent zeros are equidistant in the hyperbolic metric \eqref{hyperbolic}.
	\end{Corollary}
	
	\begin{proof}
		Suppose that $\mathcal{T}_{z^n} \cong \mathcal{T}_B$.
		Noting that $(z^n)^{\#} = z^n$ and applying Theorem \ref{Main-SI}
		we conclude that $B = \psi \circ \phi^n$
		for some $\phi, \psi \in \Aut(\D)$. If $\psi$ is a rotation then
		$B$ has one zero of order $n$. If $\psi$ is not a rotation, then the
		zeros $z_1,z_2, \ldots, z_n$ of $B$ are distinct and satisfy the equation
		\begin{equation*}
			\left(\frac{z_j - a}{1 - \overline{a} z_j}\right)^n = b
		\end{equation*}
		for some $a,b \in \D$.
		The $n$th roots of $b$ are equally spaced on a circle of radius
		$|b|$ centered at the origin (which are also equally spaced with respect to
		the hyperbolic metric). The $z_j$ are formed by applying a disk automorphism
		to these $n$th roots of $b$ and thus, by the conformal invariance of the hyperbolic metric,
		are equally spaced (in the hyperbolic metric) points on some
		circle $\Gamma$ in $\D$.

		Now assume that the zeros $z_1, z_2, \ldots, z_n$ of $B$ satisfy the
		hypothesis above.  If $z_1 = z_2 = \cdots = z_n$, then $B$ is the $n$th power
		of a disk automorphism and hence belongs to $\mathcal{O}(z^n)$.  In this case, we conclude that
		$\T_B \cong \T_{z^n}$.  In the second case, map the hyperbolic center
		of the circle $\Gamma$ to the origin
		with a disk automorphism $\psi$. The map $\psi$ will also map the circle
		$\Gamma$ to a circle $|z| = r$ having the same hyperbolic radius as $\Gamma$.
		Consequently, $\psi$ will map the zeros of $B$ to points
		$t_1,t_2, \ldots, t_n$ on $|z| = r$ which are equally spaced in the hyperbolic metric.
		By basic properties of the hyperbolic metric, these points take the
		form $t_j = \omega^j a$ where $\omega$ is a primitive $n$th root of unity
		and $a \in \D$. Putting this all together, we get that the
		zeros $z_1,z_2, \cdots, z_n$ of $B$ satisfy
		\begin{equation*}
			z_j = \psi^{-1}(w^j a)
		\end{equation*}
		and hence
		\begin{equation*}
			B = \frac{\psi^n - a}{1 - \overline{a} \psi^n} \in \mathcal{O}(z^n).
		\end{equation*}
		By Theorem \ref{Main-SI} we conclude that $\T_B \cong \T_{z^n}$. \qedhere
	\end{proof}
	
	\begin{Remark}
For any inner function $\Theta$, a well-known theorem of Frostman \cite[p.~79]{MR2261424} implies there are many $\psi \in \mbox{Aut}(\D)$ for which $B = \psi \circ \Theta$ is a Blaschke product. An application of Theorem \ref{Main-SI} shows that $\mathcal{T}_{\Theta} \cong \mathcal{T}_{B}$. It is natural to ask whether or not there are infinite Blaschke products $B$ for which $\mathcal{T}_{\Theta} \cong \mathcal{T}_{B}$ implies that $\Theta$ is a Blaschke product. Again, using Theorem \ref{Main-SI}, this can be rephrased as: for a fixed infinite Blaschke product $B$, when does $\mathcal{O}(B) \cup \mathcal{O}(B^{\#})$ contain only Blaschke products? A little exercise will show that this is true precisely when $\psi \circ B$ is a Blaschke product for every $\psi \in \mbox{Aut}(\D)$. Blaschke products satisfying this property are called \textit{indestructible} (see \cite{IBP} and the references therein). It is well-known that \textit{Frostman Blaschke products}
 i.e., those Blaschke products $B$ which satisfy
		\begin{equation*}
			\sup_{\zeta \in \dD} \sum_{n = 1}^{\infty} \frac{1 - |a_{n}|^2}{|\zeta - a_n|} < \infty,
		\end{equation*}
		where $(a_n)_{n \geq 1}$ are the zeros of $B$, repeated accordingly to multiplicity,
 are indestructible. Moreover, using a deep theorem of Hruscev and Vinogradov concerning the inner multipliers of the space of Cauchy transforms of measures on the unit circle \cite[Ch.~6]{MR2215991} along with a result from \cite{MR2321805}, one can show that $\mathcal{O}(B) \cup \mathcal{O}(B^{\#})$ contains only Frostman Blaschke products if and only if $B$ is a Frostman Blaschke product.

	\end{Remark}

\section{Proof of Theorem \ref{Main-SI}}\label{SectionProof}

	The proof of Theorem \ref{Main-SI} is somewhat lengthy
	and it is consequently broken up into a series of propositions and lemmas.
	For the sake of clarity, we deal with the implications $(\Leftarrow)$ and $(\Rightarrow)$
	in equation \eqref{eq-Implication} separately.

\subsection*{Proof of the implication $(\Leftarrow)$ in \eqref{eq-Implication}}
	This is the simpler portion of the proof and it boils down to several
	computational results.
	
	\begin{Proposition} \label{P-COV}
		If $\Theta$ is inner and $\psi \in \Aut(\D)$, then
		$\T_{\Theta} \cong \T_{\Theta \circ \psi}$.
	\end{Proposition}

	\begin{proof}
		Let
		\begin{equation}\label{eq-PD}
			\psi(z) = \eta \frac{z - a}{1 - \overline{a} z}, \quad \eta \in \dD, a \in \D,
		\end{equation}
		be a typical disk automorphism and define $U:H^2 \to H^2$ by
		\begin{equation*}
			U f = \sqrt{\psi'} (f \circ \psi).
		\end{equation*}
		One can check by the change of variables formula that $U$ is a unitary operator and
		\begin{equation*}
			U^{-1} f = U^{*} f = \sqrt{(\psi^{-1})'} (f \circ \psi^{-1}).
		\end{equation*}
		Next observe that if $f \in \mathcal{K}_{\Theta}$, then
		\begin{equation*}
			\inner{U f, (\Theta \circ \psi) h}
			= \inner{f, U^{*}((\Theta \circ \psi) h)}
			= \inner{f, \Theta \sqrt{(\psi^{-1})'} (h \circ \psi^{-1})} = 0
		\end{equation*}
		for all $h \in H^2$. Similarly, for $g \in \mathcal{K}_{\Theta \circ \psi}$ we have
		\begin{equation*}
			\inner{U^{*} g, \Theta h}
			= \inner{g, U(\Theta h)} = \inner{g, \sqrt{\psi'} (\Theta \circ \psi) (h \circ \psi)} = 0.
		\end{equation*}
		Thus $U\K_{\Theta} = \K_{\Theta \circ \psi}$ and hence
		$U$ restricts to a unitary map from $\K_{\Theta}$ onto $\K_{\Theta \circ \psi}$,
		which we also denote by $U$.
		
		If $A_{g}^{\Theta} \in \mathcal{T}_{\Theta}$, then
		observe that for $f \in \mathcal{K}_{\Theta}$ we have
		\begin{align*}
			[U A_{g}^{\Theta} f](\lambda)
			& =  [U P_{\Theta} (g f)](\lambda)\\
			& =  \sqrt{\psi'(\lambda)} [P_{\Theta}(g f)](\psi(\lambda))\\
			&=  \sqrt{\psi'(\lambda)} \inner{gf, k_{\psi(\lambda)}}   && \text{by \eqref{eq-ProjectionTheta}}\\
			& =  \sqrt{\psi'(\lambda)} \int_{\dD} g(\zeta) f(\zeta)
			\left(\frac{1 - \Theta(\psi(\lambda)) \overline{\Theta(\zeta)}}
				{1 - \overline{\zeta} \psi(\lambda)}\right) \frac{|d \zeta|}{2 \pi}
\end{align*}
Now make the change of variables $\zeta = \psi(w)$ and use the identities
$$\psi(z) = \eta \frac{z - a}{1 - \overline{a} z}, \quad \psi'(z) = \eta \frac{1 - |a|^2}{(1 - \overline{a} z)^2}$$
to show that the above is equal to
			$$\int_{\dD} \sqrt{\psi'(w)} g(\psi(w)) f(\psi(w))
			\frac{1 - \Theta(\psi(\lambda)) \overline{\Theta(\psi(w))}}
			{1 - \overline{w} \lambda} \frac{|d w|}{2 \pi} = [A_{g \circ \psi}^{\Theta \circ \psi} U f](\lambda)  .
$$
		From this we conclude that $U A_{g}^{\Theta} = A_{g \circ \psi}^{\Theta \circ \psi} U$
		whence $A \mapsto U A U^{*}$ is a spatial isomorphism
		between $\mathcal{T}_{\Theta}$ and $\mathcal{T}_{\Theta \circ \psi}$.
	\end{proof}

	The computational portion of the following proposition is originally due to Crofoot \cite{Crofoot}.
	A detailed discussion of these so-called \emph{Crofoot transforms}
	in the context of truncated Toeplitz operators can be found in \cite[Sec.~13]{Sarason}.

	\begin{Proposition}[Crofoot] \label{P-Crofoot}
		If $\Theta$ is inner, $a \in \D$, and $\phi_{a}$ denotes the M\"obius transformation
		\eqref{eq-Mobius}, then
		\begin{equation*}
			U f := \frac{\sqrt{1 - |a|^2}}{1 - \overline{a} \Theta} f
		\end{equation*}
		defines a unitary operator from $\mathcal{K}_{\Theta}$ to $\mathcal{K}_{\phi_a \circ \Theta}$.
		Moreover, $U \mathcal{T}_{\Theta} U^{*} = \mathcal{T}_{\phi_a \circ \Theta}$.
		Thus for any $\psi \in \Aut(\D)$ we have
		$\mathcal{T}_{\Theta} \cong \mathcal{T}_{\psi \circ \Theta}$.
	\end{Proposition}

	Our next goal is to establish that $\T_{\Theta} \cong \T_{\Theta^{\#}}$.  This is
	the content of Proposition \ref{P-sharp} below. We should remark that this observation
	is closely related to \cite[Cor.~1.7, Prop.~1.8]{Ber}.  The proof of Proposition \ref{P-sharp}
	requires two preliminary lemmas.
	First, recall the definitions \eqref{eq-ModelConjugation} of the
	conjugation $C$ on the model space $\mathcal{K}_{\Theta}$ and \eqref{eq-Sharp}
	of the conjugation $f \mapsto f^{\#}$.
	Now let $C^{\#}$ denote the corresponding conjugation on the model space $\mathcal{K}_{\Theta^{\#}}$.
	Finally, we define a conjugate-linear map $J$ on $\mathcal{K}_{\Theta}$ by $Jf = f^{\#}$.

	\begin{Lemma} \label{bunch-o-facts}
		For $\Theta$ inner,
		\begin{enumerate} \addtolength{\itemsep}{0.5\baselineskip}
			\item $J\mathcal{K}_{\Theta} = \mathcal{K}_{\Theta^{\#}}$.

			\item If $g \in \mathcal{K}_{\Theta^{\#}}$, then $J^{-1} g = g^{\#}$.

			\item $J C: \mathcal{K}_{\Theta} \to \mathcal{K}_{\Theta^{\#}}$ is unitary.
				Also, the following formulae hold
				\begin{equation*}
					J C = C^{\#} J, \quad (J C)^{*} = C J^{-1} = J^{-1} C^{\#}.
				\end{equation*}

			\item For all $\lambda \in \D$, we have
				\begin{equation*}
					J C K_{\lambda}^{\Theta} = C^{\#} K_{\overline{\lambda}}^{\Theta^{\#}},
					\quad JC (C K_{\lambda}^{\Theta}) = K_{\overline{\lambda}}^{\Theta^{\#}}.
				\end{equation*}
		\end{enumerate}
	\end{Lemma}
	
	\begin{proof}
		Statement (i) follows from the fact that $f \mapsto f^{\#}$ is a conjugation on $H^2$ and hence
		\begin{equation*}
			0 = \inner{f, \Theta h} = \inner{\Theta^{\#} h^{\#}, f^{\#}},
			\quad f \in \mathcal{K}_{\Theta}, h \in H^2.
		\end{equation*}
		Statement (ii) is immediate since $f \rightarrow f^{\#}$ is an involution on $H^2.$
		For (iii), it is clear that $JC$ is unitary since $J$ and $C$ are isometric and conjugate
		 linear. The remaining identities  in (iii) can be easily checked.
		For (iv) first compute $J K_{\lambda}^{\Theta} = K_{\overline{\lambda}}^{\Theta^{\#}}$ and finish by
		 using $JC = C^{\#} J$.
	\end{proof}

	\begin{Lemma}\label{L-Jsharp}
		If $A_{\phi}^{\Theta} \in \mathcal{T}_{\Theta}$, then
		$J A_{\phi}^{\Theta} J^{-1} = A_{\phi^{\#}}^{\Theta^{\#}}$.
	\end{Lemma}
	
	\begin{proof}
		For all $f, g \in \mathcal{K}^{\infty}_{\Theta^{\#}} $,
		\begin{align*}
			\inner{J A_{\phi}^{\Theta} J^{-1} f, g}
			&=  \inner{J g, A_{\phi}^{\Theta} J^{-1} f}\\
			& =  \inner{(A_{\phi}^{\Theta})^{*} J g, J^{-1} f}\\
			& =  \int_{0}^{2 \pi} \overline{\phi(e^{i \theta})}\overline{ g(e^{-i \theta})} f(e^{-i \theta}) \frac{ d \theta}{2 \pi}\\
			& =  \int_{0}^{-2 \pi} \overline{\phi(e^{-i \theta})} \overline{g(e^{i \theta})} f(e^{i \theta}) \frac{-d \theta}{2 \pi}\\
			& =  \int_{0}^{2 \pi} \overline{\phi(e^{-i \theta})} \overline{g(e^{i \theta})} f(e^{i \theta}) \frac{d \theta}{2 \pi}\\
			& =  \inner{\phi^{\#} f, g}\\
			& =  \inner{A_{\phi^{\#}}^{\Theta^{\#}} f, g}. \qedhere
		\end{align*}
	\end{proof}

	Armed now with Lemmas \ref{bunch-o-facts} and \ref{L-Jsharp} we are ready to prove the following.

	\begin{Proposition} \label{P-sharp}
		For $\Theta$ inner, $\mathcal{T}_{\Theta} \cong \mathcal{T}_{\Theta^{\#}}$.
	\end{Proposition}

	\begin{proof}
		From Lemma \ref{bunch-o-facts}, the operator
		\begin{equation} \label{eq-UisJC}
			JC: \mathcal{K}_{\Theta} \to \mathcal{K}_{\Theta^{\#}},
		 \end{equation}
		 is unitary. Furthermore, for $f, g \in \mathcal{K}_{\Theta^{\#}}$ we have
		\begin{align*}
			\inner{(J C) A_{\phi}^{\Theta} (J C)^{*} f, g}
			& =  \inner{C^{\#} J A_{\phi}^{\Theta} J^{-1} C^{\#} f, g}
				&& \mbox{(by Lemma \ref{bunch-o-facts})}\\
			& =  \inner{C^{\#} A_{\phi^{\#}}^{\Theta^{\#}} C^{\#} f, g}
				&& \mbox{(by Lemma \ref{L-Jsharp})}\\
			& =  \inner{(A_{\phi^{\#}}^{\Theta^{\#}})^{*} f, g}
				&& \mbox{(Proposition \ref{TTOareCSO})}\\
			& =  \inner{A_{\overline{\phi^{\#}}}^{\Theta^{\#}} f, g}.
		\end{align*}
		It follows that $A \mapsto (JC) A (JC)^{*}$
		is a spatial isomorphism from $\mathcal{T}_{\Theta}$ onto $\mathcal{T}_{\Theta^{\#}}$.
	\end{proof}

	Propositions \ref{P-COV}, \ref{P-Crofoot}, and \ref{P-sharp}
	yield the implication $(\Leftarrow)$ of \eqref{eq-Implication}.  This completes
	the first part of the proof of Theorem \ref{Main-SI}.

\subsection*{Technical Lemmas}
	The proof of the $(\Rightarrow)$ implication in \eqref{eq-Implication} is significantly
	more involved than the proof of $(\Leftarrow)$.  We require several additional technical lemmas
	which we present in this subsection.

	\begin{Lemma} \label{L-L-vs-CL}
		Let $\Theta$ be inner, $\Theta \not \in \Aut(\D)$, and let
		\begin{align*}
			\mathcal{L}_{\Theta} &:= \left\{\rho k_{\lambda}: \rho \in \dD, \lambda \in \D\right\}, \\
			\widetilde{\mathcal{L}}_{\Theta} &:= \left\{\rho C k_{\lambda}: \rho \in \dD, \lambda \in \D\right\}.
		\end{align*}
		For each fixed $\lambda_0 \in \D$, we have
		\begin{align}
			\dist\left(k_{\lambda_0}, \widetilde{\mathcal{L}}_{\Theta}\right) &> 0 ,  \label{eq-FirstStuff} \\
			\dist\left(C k_{\lambda_0}, \mathcal{L}_{\Theta}\right) &> 0. \label{eq-SecondStuff}
		\end{align}
	\end{Lemma}
	
	\begin{proof}
		Suppose that $\dist\left(k_{\lambda_0}, \widetilde{\mathcal{L}}_{\Theta}\right) = 0$
		holds for some $\lambda_0 \in \D$.  It follows that
		there are sequences $(\mu_n)_{n \ge 1} \subset \D$ and $(\rho_{n})_{n \geq 1} \subset \dD$ so that
		\begin{equation}\label{eq-RhoH}
			\rho_n C k_{\mu_n} \to k_{\lambda_0}
		\end{equation}
		in the norm of $H^2$.  Passing to a subsequence, we can assume that $\mu_n$
		converges to some $\mu_0 \in \D^{-}$.  There are two cases we must consider.
		\medskip
		
		\noindent\textbf{Case 1}:  If $\mu_0 \in \D$, then
		\begin{equation*}
		C k_{\mu_n} \to C k_{\mu_0}
		\end{equation*}
		in $H^2$ and hence pointwise in $\D$. This forces the sequence $\rho_n$ to converge to some $\rho_0 \in \dD$ and hence
		\begin{equation*}
		k_{\lambda_0} = \rho_0 C k_{\mu_0}.
		\end{equation*}
		However, this contradicts Lemma \ref{L-LI} from which we conclude that $\mu_0 \in \dD$.
		\medskip
		
		\noindent\textbf{Case 2}:  If $\mu_0 \in \dD$, then the sequence $\Theta(\mu_n)$ is bounded
		and hence upon passing to a subsequence we may assume that $\Theta(\mu_n) \to a$
		for some $a \in \D^{-}$.  By \eqref{eq-RhoH} it follows that
		\begin{equation} \label{eq-rescue}
			\rho_n \frac{\Theta(z) - \Theta(\mu_n)}{(z - \mu_n) \|C K_{\mu_n}\|}
			\overset{H^2}{\longrightarrow}
			\frac{1 - \overline{\Theta(\lambda_0)} \Theta(z)}{(1 - \overline{\lambda_0} z) \|K_{\lambda_0}\|}
		\end{equation}
		whence we also have pointwise convergence on $\D$.  For any fixed $z_0 \in \D$
		for which $\Theta(z_0) \not = a$ we conclude that
		\begin{equation*}
				\rho_n \frac{\Theta(z_0) - \Theta(\mu_n)}{(z_0 - \mu_n) \|C K_{\mu_n}\|}
				\to \frac{1 - \overline{\Theta(\lambda_0)} \Theta(z_0)}{(1 - \overline{\lambda_0} z_0)
				\|K_{\lambda_0}\|} \not = 0.
		\end{equation*}
		But since
		\begin{equation*}
			\frac{\Theta(z_0) - \Theta(\mu_n)}{z_0 - \mu_n} \to \frac{\Theta(z_0) - a}{z_0 - \mu_0} \not = 0,
		\end{equation*}
		it follows that $\rho_n$ converges to some  $\rho_0 \in \dD$ and $\|C K_{\mu_n}\|^{-1}$ converges
		to some finite number $M$.  Upon letting $n \to \infty$ in \eqref{eq-rescue}, we obtain
		\begin{equation*}
			\rho_0 M \frac{\Theta(z) - a}{z - \mu_0}
			= \frac{1 - \overline{\Theta(\lambda_0)} \Theta(z)}{(1 - \overline{\lambda_0} z) \|K_{\lambda_0}\|}.
		\end{equation*}
		Solving for $\Theta(z)$ in the preceding reveals that $\Theta$ is a linear fractional transformation
		-- contradicting the assumption that $\Theta \not \in \mbox{Aut}(\D)$.
		This establishes \eqref{eq-FirstStuff}.  The second inequality \eqref{eq-SecondStuff}
		follows immediately since $C$ is an involutive isometry and so
		\begin{equation*}
			\mbox{dist}(k_{\lambda_0}, \widetilde{\mathcal{L}}_{\Theta})
			 =  \mbox{dist}(C k_{\lambda_0}, C^2 \mathcal{L}_{\Theta})
			 =  \mbox{dist}(C k_{\lambda_0}, \mathcal{L}_{\Theta}).\qedhere
		\end{equation*}
	\end{proof}

	We henceforth assume that $\Theta_1$ and $\Theta_2$ are fixed inner functions, neither in $\mbox{Aut}(\D)$,
	and that $U \mathcal{T}_{\Theta_1} U^{*} = \mathcal{T}_{\Theta_2}$ for some unitary
	$U: \mathcal{K}_{\Theta_1} \to \mathcal{K}_{\Theta_2}$.
	We let $C_1,C_2$ denote the conjugations \eqref{eq-ModelConjugation} on
	$\K_{\Theta_1}$ and $\K_{\Theta_2}$, respectively.
	To simplify our notation somewhat, we set
	\begin{equation*}
		k_{\lambda} := k_{\lambda}^{\Theta_1}, \quad \widetilde{k}_{\lambda} := C_1 k_{\lambda}^{\Theta_1},
		\quad  \ell_{\lambda} := k_{\lambda}^{\Theta_2}, \quad \widetilde{\ell}_{\lambda} := C_2 k_{\lambda}^{\Theta_2}
	\end{equation*}
	for $\lambda \in \D$.

	We now exploit the fact that the rank-one operators in $\mathcal{T}_{\Theta_1}$
	are carried onto the rank-one operators in $\mathcal{T}_{\Theta_2}$ by our spatial isomorphism.
	By Lemma \ref{Sarason-tensors-kck} and Lemma \ref{L-LI}, we conclude that
	$U (k_{\lambda} \otimes \widetilde{k}_{\lambda}) U^{*}$ is either
	$\zeta \ell_{\eta} \otimes \widetilde{\ell}_{\eta}$ for some $\zeta \in \dD$ and $\eta \in \D$, or
	$\zeta' \widetilde{\ell}_{\eta'} \otimes \ell_{\eta'}$ for some $\zeta' \in \dD$ and $\eta' \in \D$.
	Upon applying Lemma \ref{L-tensors} we observe that
	\begin{equation}\label{eq-EitherOr}
		U k_{\lambda} \in \mathcal{L}_{\Theta_2} \cup \widetilde{\mathcal{L}}_{\Theta_2}.
	\end{equation}
	
	In fact, even more is true.

	\begin{Lemma} \label{C-L-vs-CL}
		Either $U \mathcal{L}_{\Theta_1} = \mathcal{L}_{\Theta_2}$
		or $U \mathcal{L}_{\Theta_1} = \widetilde{\mathcal{L}}_{\Theta_2}$.
		As a consequence, there are maps $w: \D \to \dD$ and $\phi: \D \to \D$ so that either
		\begin{equation*}
			U(k_{\lambda} \otimes \widetilde{k}_{\lambda}) = w(\lambda) \ell_{\phi(\lambda)}
			\otimes \widetilde{\ell}_{\phi(\lambda)}, \quad \forall \lambda \in \D,
		\end{equation*}
		or
		\begin{equation*}
			U(k_{\lambda} \otimes \widetilde{k}_{\lambda}) = w(\lambda) \widetilde{\ell}_{\phi(\lambda)}
			\otimes \ell_{\phi(\lambda)}, \quad \forall \lambda \in \D.
		\end{equation*}
	\end{Lemma}
	
	\begin{proof}
		Since the map $\lambda \mapsto k_{\lambda}$ is continuous
		from $\D$ to $\mathcal{K}_{\Theta_1}$, it follows that
		\begin{equation*}
			F(\lambda) := U k_{\lambda}
		\end{equation*}
		is a continuous function from $\D$ to $\mathcal{K}_{\Theta_2}$.
		Suppose that $F(\lambda_0) = \rho_0 \ell_{\eta_0} \in \mathcal{L}_{\Theta_2}$
		for some $\lambda_0, \eta_0 \in \D, \rho_0 \in \dD$.
		We now show that there is an open disk $B(\lambda_0, \delta)$ about $\lambda_0$ (of radius $\delta>0$) so that
		\begin{equation*}
			\lambda \in B(\lambda_0, \delta)
			\quad \Rightarrow \quad
			U k_{\lambda} \in \mathcal{L}_{\Theta_2}.
		\end{equation*}
		If this were not the case then by \eqref{eq-EitherOr} there exists
		sequences $\lambda_n \to \lambda_0, \eta_n \in \D, \rho_n \in \dD$ so that
		$F(\lambda_n) = \rho_n \widetilde{\ell}_{\eta_n}$.
		By the continuity of $F$ at $\lambda_0$, we see that
		$\rho_n \widetilde{\ell}_{\eta_n} \to \rho_0 \ell_{\eta_0}$, which contradicts Lemma \ref{L-L-vs-CL}.
		Since $\D$ is connected, we conclude that $U \mathcal{L}_{\Theta_1} \subset \mathcal{L}_{\Theta_2}$.

		If we now interchange the roles of $\Theta_1$ and $\Theta_2$, replacing $U$ with $U^{*}$, the
		argument above shows that $U^{*} \mathcal{L}_{\Theta_1} \subset \mathcal{L}_{\Theta_1}$.
		This means that $\mathcal{L}_{\Theta_2} \subset U \mathcal{L}_{\Theta_1}$
		and so $U \mathcal{L}_{\Theta_1} = \mathcal{L}_{\Theta_2}$.
		The same argument shows that if $F(\lambda_0) \in \widetilde{\mathcal{L}}_{\Theta_2}$,
		then $U \mathcal{L}_{\Theta_1} = \widetilde{\mathcal{L}}_{\Theta_2}$.
	\end{proof}

	\begin{Remark} \label{e-o-2}
		Now observe that it suffices to consider the case where
		$U \mathcal{L}_{\Theta_1} = \mathcal{L}_{\Theta_2}$.
		Indeed, suppose that $U \mathcal{L}_{\Theta_1} = \widetilde{\mathcal{L}}_{\Theta_2}$.
		We know from Proposition \ref{P-sharp} that
		$\mathcal{T}_{\Theta_2} \cong \mathcal{T}_{\Theta_{2}^{\#}}$ and, from
		Lemma \ref{bunch-o-facts} part (iv), the unitary $JC$ implementing this
		spatial isomorphism carries $\widetilde{\mathcal{L}}_{\Theta_2}$ onto
		$\mathcal{L}_{\Theta_{2}^{\#}}$.  By replacing $\Theta_2$ with $\Theta_{2}^{\#}$
		if necessary (which does not change $\mathcal{O}(\Theta_2) \cup \mathcal{O}(\Theta_{2}^{\#})$), we assume for the remainder of the proof that $U \mathcal{L}_{\Theta_1} = \mathcal{L}_{\Theta_2}$.
		Under this assumption it follows that
		\begin{equation} \label{e-auto-F}
			U(k_{\lambda} \otimes \widetilde{k}_{\lambda})
			= w(\lambda) \ell_{\phi(\lambda)} \otimes \widetilde{\ell}_{\phi(\lambda)}, \quad \forall \lambda \in \D,
		\end{equation}
		for some functions $w: \D \to \dD$ and $\phi: \D \to \D$.
	\end{Remark}

	\begin{Lemma} \label{L-auto-F}
		The function $\phi$ in \eqref{e-auto-F} belongs to $\Aut(\D)$.
	\end{Lemma}

	\begin{proof}
		We first prove that $\phi:\D\to\D$ is a bijection.
		Suppose that $\phi(\lambda_1) = \phi(\lambda_2)$.
		It follows from \eqref{e-auto-F} that $k_{\lambda_1} = c k_{\lambda_2}$ for some scalar $c$.
		By Lemma \ref{L-LI}, we conclude that $\lambda_1 = \lambda_2$ whence $\phi$ is injective.
		Now let $\eta \in \D$.  By Lemma \ref{C-L-vs-CL} we know that
		\begin{equation*}
			U^{*} (\ell_{\eta} \otimes \widetilde{\ell}_{\eta}) U = c k_{\lambda} \otimes \widetilde{k}_{\lambda}
		\end{equation*}
		for some $\lambda \in \D$ and some scalar $c$.   We cannot have
		\begin{equation*}
			U^{*} (\ell_{\eta} \otimes \widetilde{\ell}_{\eta}) U = c \widetilde{k}_{\lambda} \otimes k_{\lambda}
		\end{equation*}
		or else (by Lemma \eqref{L-tensors}) $U k_{\lambda} = c \widetilde{\ell}_{\eta}$
		which we are assuming is not the case. Another application of
		Lemma \ref{L-tensors} reveals that $\phi(\lambda) = \eta$ whence $\phi$ is surjective.

		To show that $\phi \in \Aut(\D)$, it suffices to prove that $\phi$ is analytic on $\D$.
		We may assume that $\Theta_2(0) \not = 0$ and $\Theta_2(w_0) = 0$ for some
		$w_0 \in \D$. If this is not the case, choose $a_1, a_2 \in \D$ ($a_1 \neq a_2$)
		so that $\Theta_2(a_1) = \Theta_2(a_2) = b$, replace $\Theta_2$ by
		$\phi_b \circ \Theta_2 \circ \phi_{-a}$, and appeal to Propositions \ref{P-COV} and \ref{P-Crofoot}.
		In particular, this means that if $L_{\eta}$ denotes the
		reproducing kernel for $\mathcal{K}_{\Theta_2}$, then
		\begin{equation*}
			L_0 = 1, \quad L_{w_0} = \frac{1}{1 - \overline{w_0} z}.
		\end{equation*}
		Let $f = U^{-1} L_0$ and $g = U^{-1} L_{w_0}$. Then for any $\lambda \in \D$ we have
		\begin{align*}
			f(\lambda)
			& =  \inner{f, K_{\lambda}}\\
			& =  \inner{U f, U K_{\lambda}}\\
			& =  \inner{1, \frac{w(\lambda) \norm{K_{\lambda}}}{\norm{L_{\phi(\lambda)}}} L_{\phi(\lambda)}}\\
			& =  \frac{\overline{w(\lambda)} \norm{K_{\lambda}}}{\norm{L_{\phi(\lambda)}}}.
		\end{align*}
		Similarly, using the formula for $f(\lambda)$ above, we get
		\begin{align*}
			g(\lambda)
			& =  \inner{g, K_{\lambda}}\\
			& =  \inner{U g, U K_{\lambda}}\\
			& =  \inner{\frac{1}{1 - \overline{w_0} z},
				\frac{w(\lambda) \norm{K_{\lambda}}}{\norm{L_{\phi(\lambda)}}} L_{\phi(\lambda)}}\\
			& =  \frac{1}{1 - \overline{w_0} \phi(\lambda)}
				 \frac{\overline{w(\lambda)} \norm{K_{\lambda}}}{\norm{L_{\phi(\lambda)}}}\\
			& =  \frac{1}{1 - \overline{w_0} \phi(\lambda)} f(\lambda).
		\end{align*}
		Since the functions $f$ and $g$ are analytic (and not identically zero) on $\D$, upon
		solving for $\phi(\lambda)$ in the preceding identity we conclude that $\phi$ is analytic on $\D$.
	\end{proof}

\subsection*{Proof of the implication $(\Rightarrow)$ in \eqref{eq-Implication}}
	We have already seen via Propositions \ref{P-COV}, \ref{P-Crofoot}, and \ref{P-sharp} that
	\begin{equation*}
		\Theta_1 \in \mathcal{O}(\Theta_2) \cup \mathcal{O}(\Theta_2^{\#})
		\quad \Rightarrow \quad \mathcal{T}_{\Theta_1} \cong \mathcal{T}_{\Theta_2}.
	\end{equation*}
We now prove the reverse implication.
	In light of Remark \ref{e-o-2} and Lemma \ref{L-auto-F} we may assume that
	\begin{equation} \label{AA}
		U k_{\lambda} = w(\lambda) \ell_{\phi(\lambda)}, \quad \forall \lambda \in \D,
	\end{equation}
	for some functions $w: \D \to \dD$ and $\phi \in \mbox{Aut}(\D)$.  Consequently
	we may appeal to Lemma \ref{L-tensors} to conclude that
	\begin{equation*}
		U \left(k_{\lambda} \otimes \widetilde{k}_{\lambda}\right) U^{*}
		= w(\lambda) \ell_{\phi(\lambda)} \otimes \widetilde{\ell}_{\phi(\lambda)}.
	\end{equation*}
	Upon taking adjoints in the preceding equation we then obtain
	\begin{equation*}
		U \left(\widetilde{k}_{\lambda} \otimes k_{\lambda}\right) U^{*}
		= \overline{w(\lambda)}\; \widetilde{\ell}_{\phi(\lambda)} \otimes \ell_{\phi(\lambda)}.
	\end{equation*}
	Lemma \ref{L-tensors} now yields
	\begin{equation} \label{BB}
		U \widetilde{k}_{\lambda} = w(\lambda) \widetilde{\ell}_{\phi(\lambda)}.
	\end{equation}

	Next we combine \eqref{AA} and \eqref{BB} to obtain
	\begin{equation*}
		|\langle\widetilde{k}_{\lambda}, k_{\lambda}\rangle|
		= |\langle\widetilde{\ell}_{\phi(\lambda)}, \ell_{\phi(\lambda)}\rangle|.
	\end{equation*}
	Noting that
	\begin{equation*}
		k_{\lambda} = \frac{K_{\lambda}}{\norm{K_{\lambda}}},
		\qquad \inner{C K_{\lambda}, K_{\lambda}}
		= \Theta'(\lambda), \qquad \|C K_{\lambda}\|
		= \|K_{\lambda}\| = \sqrt{\frac{1 - |\Theta(\lambda)|^2}{1 - |\lambda|^2}}
	\end{equation*}
	we get
	\begin{equation}\label{eq-BI}
		\frac{|\Theta_{1}'(\lambda)| (1 - |\lambda|^2)}{1 - |\Theta_1(\lambda)|^2}
		= \frac{|\Theta_{2}'(\phi(\lambda))| (1 - |\phi(\lambda)|^2)}{1 - |\Theta_2(\phi(\lambda))|^2}.
	\end{equation}
	Using the Schwarz-Pick lemma \cite[p.~2]{MR2261424} we have
	\begin{equation*}
		|\phi'(z)| = \frac{1 - |\phi(z)|^2}{1 - |z|^2}, \quad \forall z \in \D, \quad \forall \phi \in \Aut(\D),
	\end{equation*}
	whence the identity \eqref{eq-BI} becomes
	\begin{equation*}
		\frac{|\Theta_{1}'(\lambda)|}{1 - |\Theta_1(\lambda)|^2} =
		\frac{|\Theta_{2}'(\phi(\lambda))|}{1 - |\Theta_2(\phi(\lambda))|^2} |\phi'(\lambda)|.
	\end{equation*}
	Replacing $\Theta_2$ by $\Theta_2 \circ \phi$ in the preceding formula gives us
	\begin{equation} \label{final}
		\frac{|\Theta_{1}'(\lambda)|}{1 - |\Theta_1(\lambda)|^2}
		= \frac{|\Theta_{2}'(\lambda)|}{1 - |\Theta_2(\lambda)|^2}.
	\end{equation}
	Another application of the Schwarz-Pick lemma shows that \eqref{final}
	continues to hold if $\Theta_1$ is replaced by $\psi \circ \Theta_1$
	for all $\psi \in \mbox{Aut}(\D)$. It follows that we may assume that
	\begin{equation} \label{conditions-zeros}
		\Theta_1(0) = \Theta_2(0) = 0, \quad \Theta_{1}'(0) \not = 0,
		\quad \Theta_{2}'(0) \not = 0.
	\end{equation}
	If not, choose $a \in \D$ so that $\Theta_{1}'(a) \not = 0$ and $\Theta_{2}'(a) \not = 0$.
	Let $b_1 = \Theta_1(\phi_{-a}(0))$ and $b_{2} = \Theta_2(\phi_{-a}(0))$. Now replace
	$\Theta_1$ by $\phi_{b_1} \circ \Theta_1 \circ \phi_{-a}$ and $\Theta_2$
	by $\phi_{b_2} \circ \Theta_2 \circ \phi_{-a}$ and observe that
	\eqref{conditions-zeros} still holds. It is important to note that all of these simplifying assumptions on $\Theta_2$ has not altered $\mathcal{O}(\Theta_2) \cup \mathcal{O}(\Theta_{2}^{\#})$.

	The assumption \eqref{conditions-zeros} means that both $\Theta_1$
	and $\Theta_2$ are invertible near the origin.  Thus there is an
	$\varepsilon > 0$ such that $\Theta_1$ and $\Theta_2$ are injective on
	the disk $B(0, \varepsilon)$. There is also a $\delta > 0$ with
	$B(0, \delta) \subset \Theta_1(B(0, \varepsilon))$ and
	$B(0, \delta) \subset \Theta_2 (B(0, \varepsilon))$.

	Now suppose that $|z| < \delta$. Then $\Theta_{1}^{-1}([0, z])$ is a curve
	$\gamma$ in $B(0, \varepsilon)$ and
	$\Theta_2 \circ \Theta_{1}^{-1}([0, z]) = \Theta_2(\gamma)$ is a
	curve $\Gamma$ in $B(0, \delta)$ going from $0$ to
	$\beta := \Theta_2 \circ \Theta_{1}^{-1}(z)$. From our discussion in the previous
	paragraph along with the change of variables formula and \eqref{final} we get
	\begin{equation*}
		\int_{\gamma} \frac{|\Theta_{1}'(t)|}{1 - |\Theta_1(t)|^2} dt
		= \int_{\gamma} \frac{|\Theta_{2}'(t)|}{1 - |\Theta_2(t)|^2} dt = \int_{\Gamma} \frac{|d w|}{1 - |w|^2}.
	\end{equation*}
	Thus $\rho(0, z) \geq \rho(0, \beta)$ whence, by \eqref{log-identity},
	$|z| \geq |\beta|$ and so
	\begin{equation*}
		|\Theta_2 \circ \Theta_{1}^{-1}(z)| \leq |z|
	\end{equation*}
	for small $|z|$. A similar argument also shows that
	$|\Theta_1 \circ \Theta_{2}^{-1}| \leq |z|$ for small $|z|$.
	Putting this all together we find that
	\begin{equation*}
		|z| = |\Theta_{2} \circ \Theta_{1}^{-1}(z)|, \quad \forall |z| < \delta
	\end{equation*}
	and hence there is a  $\zeta \in \dD$ such that
	\begin{equation*}
		\Theta_2 \circ \Theta_{1}^{-1}(z) = \zeta z, \quad \forall |z| < \delta.
	\end{equation*}
	Replacing $z$ by $\Theta_1(z)$ for $|z|$ small, we have $\Theta_2(z) = \zeta \Theta_1(z)$
	and so $\Theta_2 = \zeta \Theta_1$ on $\D$. Thus $\Theta_1 \in \mathcal{O}(\Theta_2)$ as desired.
	This completes the proof of Theorem \ref{Main-SI}.
	\qed

\section{Unitary equivalence to a truncated Toeplitz operator}\label{SectionUE}

	In this section we attempt to describe those
	classes of Hilbert space operators which are UETTO
	(unitarily equivalent to a truncated Toeplitz operator).
	This question is more subtle that it might at first appear.
	For instance, the Volterra integration operator, being
	the Cayley transform of the compressed shift $A_z$
	on a certain model space, is UETTO \cite{MR0192355} (see also \cite[p.~41]{MR1892647}).
	While the general question appears quite difficult, we are able
	to obtain concrete results in a few specific cases.


	\begin{Theorem}\label{TheoremRankOneUETTO}
		Every rank one operator is UETTO.
	\end{Theorem}
	
	\begin{proof}
		Let $T = u \otimes v$ be a rank one operator on an $n$-dimensional Hilbert space.
		Without loss of generality, suppose that $2 \leq n \leq \infty$, $\norm{u} = \norm{v} = 1$ and
		\begin{equation*}
			0 \leq \inner{u,v} \leq 1.
		\end{equation*}
		We claim that there exists a Blaschke product $\Theta$ of order $n$ (i.e., having
		$n$ zeros, counting according to multiplicity) and an appropriate
		$\lambda$ so that $u \otimes v$ is unitarily equivalent to a multiple of
		$k_{\lambda} \otimes C k_{\lambda}$.
		By Lemmas \ref{L-tensors} and \ref{Sarason-tensors-kck} it suffices to exhibit $\Theta$ and $\lambda$ so that
		\begin{equation*}
			\inner{u,v} =  \inner{k_{\lambda},
			 C k_{\lambda}}.
		\end{equation*}
		There are three cases to consider:
		 \begin{enumerate}\addtolength{\itemsep}{0.5\baselineskip}
			\item Suppose that $\inner{u,v} = 0$.  In this case let $\Theta$ be
				a Blaschke product of order $n$ having
				a repeated root at $\lambda = 0$.  Then
				\begin{equation*}
					 \inner{k_0 ,
					C k_0}  = \inner{1, \frac{\Theta}{z} }
					=  \overline{\Theta'(0)} = 0 = \inner{u,v}
				\end{equation*}
				as desired.

		\item Suppose that $\inner{u,v} = 1$.  Since $u$ and $v$ are unit vectors, it follows
			that $u = v$.  In this case, let $\Theta$ be a Blaschke product of order $n$ having an ADC
			at $\lambda = 1$ and satisfying $\Theta(1) = 1$ in the non-tangential limiting sense.
			A short computation shows that $C k_1 = k_1$ whence
			\begin{equation*}
				 \inner{k_1,
				C k_1}  = 1 = \inner{u,v}
			\end{equation*}
			as desired.

		\item Suppose that $0 < \inner{u,v} <1$.  In this case, let $\Theta$ be a Blaschke product of order $n$
			with a simple root at $\lambda = 0$ and having its remaining roots $\lambda_i$ being
			strictly positive.  In this case
			\begin{equation*}
				\inner{k_0, C k_0}
				= \Theta'(0) = \prod_{i=1}^n \lambda_i.
			\end{equation*}
			By selecting the zeros $\lambda_i$ appropriately, the preceding can be made to equal $\inner{u,v}$ as
			was required.\qedhere
		\end{enumerate}
	\end{proof}

	\begin{Theorem} \label{T-2x2}
		Every $2 \times 2$ matrix is UETTO.  In fact, if $T$ is a given $2\times 2$
		matrix and $\Theta$ is a Blaschke product of order $2$, then $\T_{\Theta}$
		contains an operator unitarily equivalent to $T$.
	\end{Theorem}

	\begin{proof}
		Let $T$ be a given $2 \times 2$ matrix and let $\Theta$ be a Blaschke
		product of order $2$.  Using the fact that a $2 \times 2$ matrix is
		unitarily equivalent to a complex symmetric matrix
		(see \cite[Cor.~3.3]{Chevrot}, \cite[Ex.~6]{G-P}, or \cite[Cor.~3]{Tener}),
		we may restrict our attention to
		the case where $T$ is complex symmetric:  $T = T^t$.
		Now observe that the subspace of $S_2(\C) \subset M_2(\C)$ consisting of all
		$2 \times 2$ complex symmetric matrices
		has dimension $3$.  Next note that part (iv) of
		Lemma \ref{Sarason-tensors-kck}  asserts that
		$\dim \T_{\Theta} = 3$ as well.
		If $\beta$ is a $C$-real orthonormal basis for $\mathcal{K}_{\Theta}$
		(see \cite[Lem.~2.6]{G} for details), then the map
		$\Phi: \T_{\Theta} \to S_2(\C)$ defined by
		$\Phi(A) = [A]_{\beta}$  is clearly injective whence its image
		contains $T$ \cite[Lem.~2.7]{G}.
	\end{proof}
	
	\begin{Corollary}
		If $\Theta_1$ and $\Theta_2$ are Blaschke products of order 2, then
		$\T_{\Theta_1} \cong \T_{\Theta_2}$.
	\end{Corollary}
	
	\begin{proof}
		The proof of Theorem \ref{T-2x2} provides a recipe for constructing spatial isomorphisms
		$\Phi_1:\T_{\Theta_1} \to S_2(\C)$ and $\Phi_2:\T_{\Theta_2} \to S_2(\C)$.
		It follows that $\Phi_2 \circ \Phi_1: \T_{\Theta_1} \to \T_{\Theta_2}$ is a spatial isomorphism.
	\end{proof}

	\begin{Theorem}
		If $N$ is an $n \times n$ normal matrix and $\Theta$ is a Blaschke product of order $n$, then
		$N$ is unitarily equivalent to an operator in $\T_{\Theta}$.
	\end{Theorem}
	
	\begin{proof}
		By the Spectral Theorem, we know that $N$ is unitarily equivalent to the diagonal matrix
		$\operatorname{diag}(\lambda_1, \lambda_2,\ldots, \lambda_n)$ where $\lambda_1,\lambda_2,\ldots,\lambda_n$
		denote the eigenvalues of $N$, repeated according to their multiplicity.
		Select a Clark unitary operator $U = U_{\alpha}$ (see \eqref{Clark-U-defn})
		and note from Theorem \ref{T-Sarason-commutes} that $U \in \T_{\Theta}$ as
		is $p(U)$ for any polynomial $p(z)$. Also note that the eigenvalues
		$\zeta_1,\zeta_2,\ldots,\zeta_n$ of $U$ have multiplicity one \cite[Thm.~3.2]{MR0301534}
		(see also \cite[Thm.~8.2]{G}).  Thus, there exists a polynomial
		$p(z)$ such that $p(\zeta_i) = \lambda_i$ for $i = 1,2,\ldots,n$. It follows that $p(U)$ is unitarily
		equivalent to $\operatorname{diag}(\lambda_1, \lambda_2,\ldots, \lambda_n)$ and hence to $N$ itself.
	\end{proof}
	
	If we are willing to sacrifice the arbitrary selection of $\Theta$, then
	the preceding can be generalized to the infinite-dimensional setting.
	To do so, we require some preliminary remarks on multiplication operators.
	For a compactly supported Borel measure $\mu$ on $\C$, we have the associated algebra
	\begin{equation}\label{eq-MultiplicationAlgebra}
		\mathcal{M}_{\mu} := \{ M_{\phi} \in \B(L^{2}(\mu)): \phi \in L^{\infty}(\mu)\}
	\end{equation}
	of multiplication operators on $L^2(\mu)$.  For each such measure we define
	the ordered pair $\kappa(\mu) = (\epsilon, n)$ where
	\begin{equation*}
		\epsilon =
		\begin{cases}
			0 & \text{if $\mu$ is purely atomic},\\
			1 & \text{otherwise},
		\end{cases}
	\end{equation*}
	and $0 \leq n \leq \infty$ denotes the number of atoms of $\mu$.	
	In terms of the function $\kappa$, the following theorem of
	Halmos and von Neumann \cite{MR0006617} (see also
	\cite[Thm.~7.51.7]{MR1721402}) describes when the algebras
	\eqref{eq-MultiplicationAlgebra} are spatially isomorphic.

	\begin{Theorem}[Halmos and  von Neumann]\label{Halmos-vN}
		For two compactly supported Borel measures $\mu_1, \mu_2$ on $\C$,
		the algebras $\mathcal{M}_{\mu_1}$ and $\mathcal{M}_{\mu_2}$
		are spatially isomorphic if and only if $\kappa(\mu_1) = \kappa(\mu_2)$.
	\end{Theorem}

	\begin{Theorem} \label{normal-UETTO}
		Every normal operator on a separable Hilbert space is UETTO.
	\end{Theorem}
	
	\begin{proof}
		If $N$ is a normal operator on a separable Hilbert space, then
		the spectral theorem asserts that $N$ is unitarily
		equivalent to $M_{\phi}: L^2(\mu) \to L^2(\mu)$ for some compactly
		supported Borel measure $\mu$ on $\C$ and some $\phi \in L^{\infty}(\mu)$.
		Let $\eta$ be a \textit{singular} probability measure on $\dD$ for
		which $\kappa(\mu) = \kappa(\eta)$.  By Theorem \ref{Halmos-vN},
		$M_{\phi}:  L^{2}(\mu) \to L^2(\mu)$ is unitarily equivalent to
		$M_{\psi}: L^{2}(\eta) \to L^2(\eta)$, for some $\psi \in L^{\infty}(\eta)$.

		By Proposition \ref{converse-Clark}, $\eta$ is a Clark measure for some
		Clark unitary operator $U_1$ on $\K_{\Theta}$ for some inner $\Theta$.
		Again by Proposition \ref{converse-Clark}, $U_1$ is unitarily equivalent to $(M_z, L^{2}(\eta))$.
		Moreover, by Theorem \ref{T-Sarason-commutes}, we also get that $U_1$ as
		well as $\psi(U_1)$ belong to  $\mathcal{T}_{\Theta}$. Finally, note that
		\begin{equation*}
			\psi(U_1) \cong(M_{\psi}, L^2(\eta)) \cong(M_{\phi}, L^2(\mu)) \cong N.
		\end{equation*}
In the previous line we use $\cong$ to denote unitary equivalence of two operators.
	\end{proof}

	\begin{Theorem} \label{inflate}
		For $k \in \mathbb{N} \cup \{\infty\}$, the $k$-fold inflation of a finite
		Toeplitz matrix is UETTO.
	\end{Theorem}

\begin{proof}
Suppose that $n \in \N$ and $A_{\psi} \in \mathcal{T}_{z^n}$, where
\begin{equation} \label{psi-inflate}
\psi(\zeta) = \sum_{m = - n + 1}^{n - 1} a_{m} \zeta^{m}
\end{equation}
is a trigonometric polynomial. In particular, the matrix of $A_{\psi}$ relative to the usual monomial basis
		$\{1, z, \ldots, z^{n - 1}\}$ for $\mathcal{K}_{z^n}$ is a Toeplitz matrix and every finite Toeplitz
		matrix arises in this manner.

For $k \in \N \cup \{\infty\}$ let $A_{\psi} \otimes I$ denote the $k$-fold inflation of $A_{\psi}$, where $I$ is the identity matrix on some $k$-dimensional Hilbert space. We will now show that $A_{\psi} \otimes I$ is UETTO. To do this let $B$ be a Blaschke product of order $k$ (Note that $k$ can be infinite). If $T_{B}$ denotes the usual Toeplitz operator on $H^2$ with symbol $B$, then
$$T_{B} (B^j \mathcal{K}_{B}) = B^{j + 1} \mathcal{K}_{B}, \quad j = 0, 1, 2, \ldots.$$ Since
$$
H^2 = \bigoplus_{j = 0}^{\infty} B^{j} \mathcal{K}_{B},
$$
we see that $T_{B}$ is unitarily equivalent to a shift of multiplicity $k$, i.e., $T_{B} \cong T_{z} \otimes I$ (This is a standard fact from operator theory \cite[p.~111]{MR1721402}). In a similar way, one shows that
$$T_{B^{m}} \cong T_{z^{m}} \otimes I, \quad m \in \mathbb{Z},$$ and so, from \eqref{psi-inflate},
$$T_{\psi(B)} \cong T_{\psi} \otimes I.$$

A short exercise using the fact that $\mathcal{K}_{B} = (B H^2)^{\perp}$ will show that
$$\mathcal{K}_{B^{n}} = \bigoplus_{j = 0}^{n - 1} B^{j} \mathcal{K}_{B}.$$ Combine this with the above discussion to show that $A_{\psi(B)}: \mathcal{K}_{B^n} \to \mathcal{K}_{B^n}$ (which is the compression of $T_{\psi(B)}$ to $\mathcal{K}_{B^n}$) is unitarily equivalent to $A_{\psi} \otimes I$.
\end{proof}


	We conclude this section with several open questions. The first two are motivated by Theorem \ref{inflate}.

	\begin{Question}
		For which truncated Toeplitz operators $A_{\phi}^{\Theta}$ and for which $k \in \mathbb{N} \cup \{\infty\}$
		is the $k$-fold inflation of $A_{\phi}^{\Theta}$ UETTO?
	\end{Question}
	
	\begin{Question}
		When is the direct sum of truncated Toeplitz operators UETTO?
	\end{Question}

	It is known that every truncated Toeplitz operator is a complex symmetric operator
	(see Definition \ref{def-CSO} and Proposition \ref{TTOareCSO}).  Moreover,
	so is the Volterra integration operator, every $2 \times 2$ matrix, and every normal
	operator  \cite{G,G-P}.  In light of the results obtained in this section, it is natural to ask the
	following:

	\begin{Question}
		Which complex symmetric operators are UETTO?
	\end{Question}

\section{Similarity to a truncated Toeplitz operator}\label{SectionSimilar}

	It was asked in \cite{MR1723841} whether or not the inverse Jordan problem
	can be solved in the class of Toeplitz matrices. That is to say, given any Jordan canonical
	form, can one find a Toeplitz matrix that is similar to this form?
	A negative answer to this question was subsequently provided by G.~Heinig \cite{MR1839449}.
	On the other hand, it turns out that the inverse Jordan structure problem is always
	solvable in the class of truncated Toeplitz operators. In fact, we get a bit more.

	\begin{Theorem} \label{Jordan-TTO}
		Every operator on a finite dimensional space is similar
		to a co-analytic truncated Toeplitz operator.
	\end{Theorem}

	\begin{proof}
		Recalling the notation \eqref{eq-Mobius}, for a finite
		Blaschke product $\Theta$, we write
		\begin{equation}\label{eq-AsIn}
			\Theta = \phi_{z_1}^{d_1} \phi_{z_2}^{d_2}\cdots \phi_{z_r}^{d_r},
		\end{equation}
		where $z_1, z_2,\ldots ,z_r$ are the distinct  zeros of $\Theta$, and
		$d := d_1 + d_2 + \cdots +d_r $ is the order of $\Theta$.  Let
		\begin{equation*}
			\mathcal{Q}  := \{A_{\overline{\psi}} \in \mathcal{T}_{\Theta}: \psi \in H^{\infty}\}
		\end{equation*}
		denote the algebra of co-analytic truncated Toeplitz operators on $\mathcal{K}_{\Theta}$.
		Note that $\mathcal{Q}$ is the set of $A_{\overline{p}}$ where $p$
		is a polynomial of degree at most $d$.

		For $1 \leq i \leq r,$ let $P_i$ be the Riesz idempotent
		corresponding to the eigenvalue $\overline{z}_i$ of $A_{\overline{z}}$
		and note that $P_i \in \mathcal{Q} $ and
		$\ran P_i = \ker (A_{\overline{z}} - \overline{z_i} I )^{d_i}$ \cite[p.~569]{Dunford}.
		From here it is easy to see that
		\begin{equation}\label{eq-Ranges}
			\ran P_i = \mathcal{K}_{\phi_{z_i}^{d_i}}
		\end{equation}
		and that an orthonormal basis for this subspace is
		\begin{equation*}
			\{k_{z_i} \phi_{z_i}^{j-1} : 1\leq j\leq d_i\}.
		\end{equation*}

		Relative to the basis above, the restriction of
		$A_{\overline{\phi}_{z_i}}$  to $\mathcal{K}_{\phi_{z_i}^{d_i}}$ has a matrix which is a
		$d_i \times d_i$ Jordan block. Thus the algebra
		\begin{equation*}
			\mathcal{Q}_i := \mathcal{Q} | \mathcal{K}_{\phi_{z_i}^{d_i}}
		\end{equation*}
		is spatially isomorphic to the algebra of $d_i \times d_i$
		upper triangular Toeplitz matrices.

		Since
		\begin{equation*}
			\mathcal{K}_{\Theta}
			= \mathcal{K}_{\phi_{z_1}^{d_1}} \oplus \mathcal{K}_{\phi_{z_2}^{d_2}}
			\oplus \cdots \oplus \mathcal{K}_{\phi_{z_r}^{d_r}},
		\end{equation*}
		is a (non-orthogonal) direct sum of vector spaces, we see from \eqref{eq-Ranges} that
		\begin{equation*}
			\mathcal{Q} = \mathcal{Q}_1 \oplus \mathcal{Q}_2
			\oplus  \cdots \oplus \mathcal{Q}_r,
		\end{equation*}
		is a (non-orthogonal) direct sum of algebras. It is now clear that given
		a Jordan canonical form, we can find a co-analytic truncated Toeplitz
		operator with that form.  The number of blocks in the form is the number of
		distinct zeros of $\Theta$ and the size of each block determines
		the multiplicity of each given zero.
	\end{proof}

	The proof of Theorem \ref{Jordan-TTO} also proves the following corollary:
	\begin{Corollary}
		If $\Theta$ is a finite Blaschke product, $\mathcal{Q}$, the
		co-analytic truncated operators on $\mathcal{K}_{\Theta}$,
		is spatially similar to $\mathcal{Q}^{*} := \{A^{*}: A \in \mathcal{Q}\}$, the analytic
		truncated Toeplitz operators on $\mathcal{K}_{\Theta}$.
	\end{Corollary}

\begin{proof}
Observe that for each $k$, $\mathcal{Q}_k$ and $(\mathcal{Q}_k)^{*}$ are spatially isomorphic.
\end{proof}

		Theorem \ref{normal-UETTO} asserts that for a fixed inner function
		$\Theta$, $\mathcal{T}_{\Theta}$ contains many normal operators. However, they are not among the analytic (or co-analytic)
truncated Toeplitz operators except in trivial cases.

\begin{Proposition}
If $\Theta$ is inner and $A_{\phi} \in \mathcal{T}_{\Theta}$ is normal and not a multiple of the identity operator, then $\phi \not \in H^2 \cup \overline{H^2}$.
\end{Proposition}

\begin{proof}
Suppose that $\phi \in H^2$ and $A_{\phi} \in \mathcal{T}_{\Theta}$ is normal. Since $A_{\phi} = A_{P_{\Theta} \phi}$ \cite[Thm.~3.1]{Sarason}, we can assume that $\phi \in \mathcal{K}_{\Theta}$. Furthermore, if $K_{0} = 1 - \overline{\Theta(0)} \Theta$ is the reproducing kernel for $\mathcal{K}_{\Theta}$ at the origin, we have
$$A_{K_0} f = P_{\Theta}(f - f \overline{\Theta(0)} \Theta) = f, \quad f \in \mathcal{K}_{\Theta},$$ and so $A_{K_0} = I$ (this identity was observed in \cite[p.~499]{Sarason}). Since $A_{\phi}$ is normal if and only if $A_{\phi} - a I = A_{\phi - a K_0}$ is normal, we can set $a = \phi(0)/\norm{K_0}^2$ to assume that $A_{\phi}$ is normal with
 $$\phi \in \mathcal{K}_{\Theta} \quad \mbox{and} \quad \phi(0) = 0.$$
  This means that $\phi = z g$ for some $g \in H^2$, and, since $S^{*} \phi = (\phi - \phi(0))/z \in \mathcal{K}_{\Theta}$, we see that $g \in \mathcal{K}_{\Theta}$.

To show that $A_{\phi}$ cannot be normal, we will prove the inequality
$$\norm{A_{\phi}^{*} K_0} < \norm{A_{\phi} K_{0}}.$$
Observe that
$$A_{\phi} K_0 = P_{\Theta}(\phi - \overline{\Theta(0)} \Theta \phi) = \phi$$ since $\phi \in \mathcal{K}_{\Theta}$.
Now notice that
\begin{align*}
A_{\phi}^{*} K_0 & = P_{\Theta}(\overline{\phi} - \overline{\Theta(0)} \overline{\phi} \Theta)\\
& = 0 - \overline{\Theta(0)} P_{\Theta}(\overline{(z g)} \Theta)\\
& = -\overline{\Theta(0)} P_{\Theta} (C g) && \mbox{($C g = \overline{z g} \Theta$)}\\
& = -\overline{\Theta(0)} C g.
\end{align*}
Finally note that
\begin{align*}
\norm{A_{\phi}^{*} K_{0}} & = |\Theta(0)| \norm{C g}\\
& = |\Theta(0)| \norm{g} && \mbox{($C$ is isometric)}\\
& = |\Theta(0)| \norm{z g}\\
& = |\Theta(0)| \norm{\phi}\\
& < \norm{\phi} && \mbox{(since $|\Theta(0)| < 1$)}\\
& = \norm{A_{\phi} K_{0}}. \qedhere
\end{align*}
\end{proof}

%

\bibliography{SITTO}

\end{document}